# Oscillations of empirical distribution functions under dependence[*]

## Wei Biao Wu[1]

*University of Chicago*

**Abstract:** We obtain an almost sure bound for oscillation rates of empirical distribution functions for stationary causal processes. For short-range dependent processes, the oscillation rate is shown to be optimal in the sense that it is as sharp as the one obtained under independence. The dependence conditions are expressed in terms of physical dependence measures which are directly related to the data-generating mechanism of the underlying processes and thus are easy to work with.

## 1. Introduction

Let $\varepsilon'_0, \varepsilon_i, i \in \mathbb{Z}$, be independent and identically distributed (iid) random variables on the same probability space $(\Omega, \mathcal{A}, \mathbb{P})$. For $k \in \mathbb{Z}$ let

$$(1.1) \qquad X_k = g(\ldots, \varepsilon_{k-1}, \varepsilon_k),$$

where $g$ is a measurable function such that $X_k$ is a well-defined random variable. Then $\{X_k\}_{k \in \mathbb{Z}}$ forms a stationary sequence. The framework (1.1) is very general. See [12, 15, 19, 23] among others. The process $(X_k)$ is causal or non-anticipative in the sense that $X_k$ does not depend on future innovations $\varepsilon_{k+1}, \varepsilon_{k+2}, \ldots$. Causality is a reasonable assumption in practice. The Wiener-Rosenblatt conjecture says that, for every stationary and ergodic process $X_k$, there exists a measurable function $g$ and iid innovations $\varepsilon_i$ such that the distributional equality $(X_k)_{k \in \mathbb{Z}} =_{\mathcal{D}} (g(\ldots, \varepsilon_{k-1}, \varepsilon_k))_{k \in \mathbb{Z}}$ holds; see [13, 20]. For an overview of the Wiener-Rosenblatt conjecture see [9].

Let $F$ be the cumulative distribution function of $X_k$. Assume throughout the paper that $F$ has a square integrable density $f$ with square integrable derivative $f'$. In this paper we are interested in the oscillatory behavior of empirical distribution function

$$(1.2) \qquad F_n(x) = \frac{1}{n} \sum_{i=1}^n \mathbf{1}_{X_i \leq x}, \quad x \in \mathbb{R}.$$

In particular, we shall obtain an almost sure bound for the modulus of continuity for the function $G_n(x) = \sqrt{n}[F_n(x) - F(x)]$:

$$(1.3) \qquad \Delta_n(b) = \sup_{|x-y| \leq b} |G_n(x) - G_n(y)|,$$

---

[1]Department of Statistics, University of Chicago, 5734 S. University Ave, Chicago, IL 60637, e-mail: wbwu@galton.uchicago.edu

[*]The work is supported in part by NSF Grant DMS-0478704.

where $b = b_n$ is a sequence of positive numbers satisfying

(1.4) $$b_n \to 0 \text{ and } nb_n \to \infty.$$

Under the assumption that $X_i$ are iid, there exists a huge literature on the asymptotic behavior of $\Delta_n(\cdot)$; see Chapter 14 in [14] and the references cited therein. A powerful tool to deal with the empirical distribution $F_n$ is strong approximation [1]. In comparison, the behavior of $\Delta_n(\cdot)$ has been much less studied under dependence. In this paper we shall implement the new dependence measures proposed in [23] and obtain an almost sure bound for $\Delta_n$. For a recent account of empirical processes for dependent random variables see the monograph edited by Dehling et al [4].

The rest of the paper is structured as follows. Main results on $\Delta_n(b)$ are presented in Section 2 and proved in Section 3. Section 4 contains comparisons with results obtained under independence. Some open problems are also posed in Section 4.

## 2. Main results

We first introduce some notation. For a random variable $Z$ write $Z \in \mathcal{L}^p$ $(p > 0)$ if $\|Z\|_p := (\mathbb{E}|Z|^p)^{1/p} < \infty$. Write $\|\cdot\| = \|\cdot\|_2$. Let for $k \in \mathbb{Z}$,

$$\xi_k = (\ldots, \varepsilon_{k-1}, \varepsilon_k)$$

and for $k \geq 0$ let $\xi_k^*$ be a coupled process of $\xi_k$ with $\varepsilon_0$ replaced by $\varepsilon_0'$, i.e.,

$$\xi_k^* = (\xi_{-1}, \varepsilon_0', \varepsilon_1, \ldots, \varepsilon_k).$$

We shall write $X_k^* = g(\xi_k^*)$. Let $k \geq 1$ and define the conditional cumulative distribution function $F_k(\cdot|\xi_0)$ by

(2.1) $$F_k(x|\xi_0) = \mathbb{P}(X_k \leq x|\xi_0).$$

Note that $(X_i, \xi_i)$ is stationary. Then for almost every $\xi \in \cdots \times \mathbb{R} \times \mathbb{R}$ with respect to $\mathbb{P}$ we have

(2.2) $$F_k(x|\xi_0 = \xi) = \mathbb{P}(X_k \leq x|\xi_0 = \xi) = \mathbb{P}(X_{k+i} \leq x|\xi_i = \xi)$$

for all $i \in \mathbb{Z}$. In other words, $F_k(\cdot|\xi)$ is the cumulative distribution function of the random variable $g(\xi, \varepsilon_{i+1}, \ldots, \varepsilon_{i+k})$. Assume that for all $k \geq 1$ and almost every $\xi \in \cdots \times \mathbb{R} \times \mathbb{R}$ with respect to $\mathbb{P}$ that $F_k(x|\xi_0 = \xi)$ has a derivative $f_k(x|\xi_0 = \xi)$, which is the conditional density of $X_k$ at $x$ given $\xi_0 = \xi$. By (2.2), for any $i \in \mathbb{Z}$, $f_k(x|\xi_i)$ is the conditional density of $X_{k+i}$ at $x$ given $\xi_i$. Let the conditional characteristic function

(2.3) $$\varphi_k(\theta|\xi_0 = \xi) = \mathbb{E}(e^{\sqrt{-1}\theta X_k}|\xi_0 = \xi) = \int_\mathbb{R} e^{\sqrt{-1}\theta t} f_k(t|\xi_0 = \xi)dt,$$

where $\sqrt{-1}$ is the imaginary unit. Our dependence condition is expressed in terms of the $\mathcal{L}^2$ norm $\|\varphi_k(\theta|\xi_0) - \varphi_k(\theta|\xi_0^*)\|$.

**Theorem 2.1.** *Assume that $b_n \to 0$, $\log n = O(nb_n)$ and that there exists a positive constant $c_0$ for which*

(2.4) $$\sup_x f_1(x|\xi_0) \leq c_0$$

*holds almost surely. Further assume that*

$$\sum_{k=1}^{\infty} \left[ \int_{\mathbb{R}} (1+\theta^2) \|\varphi_k(\theta|\xi_0) - \varphi_k(\theta|\xi_0^*)\|^2 d\theta \right]^{1/2} < \infty. \tag{2.5}$$

*Let*

$$\iota(n) = (\log n)^{1/2} \log \log n. \tag{2.6}$$

*Then*

$$\Delta_n(b_n) = O_{\text{a.s.}}(\sqrt{b_n \log n}) + o_{\text{a.s.}}[b_n \iota(n)]. \tag{2.7}$$

Roughly speaking, (2.5) is a short-range dependence condition. Recall that $f_k(x|\xi_0)$ is the conditional (predictive) density of $X_k$ at $x$ given $\xi_0$ and $\varphi_k(\theta|\xi_0)$ is the conditional characteristic function. So $\varphi_k(\theta|\xi_0) - \varphi_k(\theta|\xi_0^*)$ measures the degree of dependence of $\varphi_k(\cdot|\xi_0)$ on $\varepsilon_0$. Hence the summand in (2.5) quantifies a distance between the conditional distributions $[X_k|\xi_0]$ and $[X_k^*|\xi_0^*]$ and (2.5) means that the cumulative contribution of $\varepsilon_0$ in predicting future values is finite.

For the two terms in the bound (2.7), the first one $O_{\text{a.s.}}(\sqrt{b_n \log n})$ has the same order of magnitude as the one that one can obtain under independence. See Chapter 14 in [14] and Section 4. The second term $o_{\text{a.s.}}[b_n \iota(n)]$ is due to the dependence of the process $(X_k)$. Clearly, if $b_n^{1/2}(\log \log n) = o(1)$, then the first term dominates the bound in (2.7). The latter condition holds under mild conditions on $b_n$, for example, if $b_n = O(n^{-\eta})$ for some $\eta > 0$.

Let $k \geq 1$. Observe that $\varphi_k(\theta|\xi_0) = \mathbb{E}[\varphi_1(\theta|\xi_{k-1})|\xi_0]$ and

$$\mathbb{E}[\varphi_1(\theta|\xi_{k-1})|\xi_{-1}] = \mathbb{E}[\varphi_1(\theta|\xi_{k-1}^*)|\xi_{-1}] = \mathbb{E}[\varphi_1(\theta|\xi_{k-1}^*)|\xi_0]. \tag{2.8}$$

To see (2.8), write $h(\xi_{k-1}) = \varphi_1(\theta|\xi_{k-1})$. Note that $\varepsilon_i, \varepsilon_0'$, $i \in \mathbb{Z}$, are iid and $k - 1 \geq 0$. Then we have $\mathbb{E}[h(\xi_{k-1})|\xi_{-1}] = \mathbb{E}[h(\xi_{k-1}^*)|\xi_{-1}]$ since $\xi_{k-1}^*$ is a coupled version of $\xi_{k-1}$ with $\varepsilon_0$ replaced by $\varepsilon_0'$. On the other hand, we have $\mathbb{E}[h(\xi_{k-1}^*)|\xi_{-1}] = \mathbb{E}[h(\xi_{k-1}^*)|\xi_0]$ since $\varepsilon_0$ is independent of $\xi_{k-1}^*$. So (2.8) follows. Define the projection operator $\mathcal{P}_k$ by

$$\mathcal{P}_k Z = \mathbb{E}(Z|\xi_k) - \mathbb{E}(Z|\xi_{k-1}), \quad Z \in \mathcal{L}^1.$$

By the Jensen and the triangle inequalities,

$$\begin{aligned}
\|\varphi_k(\theta|\xi_0) - \varphi_k(\theta|\xi_0^*)\| &\leq \|\varphi_k(\theta|\xi_0) - \mathbb{E}[\varphi_k(\theta|\xi_0)|\xi_{-1}]\| \\
&\quad + \|\mathbb{E}[\varphi_k(\theta|\xi_0)|\xi_{-1}] - \varphi_k(\theta|\xi_0^*)\| \\
&= 2\|\mathcal{P}_0 \varphi_1(\theta|\xi_{k-1})\| \\
&\leq 2\|\varphi_1(\theta|\xi_{k-1}) - \varphi_1(\theta|\xi_{k-1}^*)\|.
\end{aligned}$$

Then a sufficient condition for (2.5) is

$$\sum_{k=0}^{\infty} \left[ \int_{\mathbb{R}} (1+\theta^2) \|\varphi_1(\theta|\xi_k) - \varphi_1(\theta|\xi_k^*)\|^2 d\theta \right]^{1/2} < \infty. \tag{2.9}$$

In certain applications it is easier to work with (2.9). In Theorem 2.2 below we show that (2.9) holds for processes $(X_k)$ with the structure

$$X_k = \varepsilon_k + Y_{k-1}, \tag{2.10}$$

where $Y_{k-1}$ is $\xi_{k-1} = (\ldots, \varepsilon_{k-2}, \varepsilon_{k-1})$ measurable. It is also a large class. The widely used linear process $X_k = \sum_{i=0}^{\infty} a_i \varepsilon_{k-i}$ is of the form (2.10). Nonlinear processes of the form

$$X_k = m(X_{k-1}) + \varepsilon_k \tag{2.11}$$

also fall within the framework of (2.10) if (2.11) has a stationary solution. A prominent example of (2.11) is the threshold autoregressive model [19]

$$X_k = a \max(X_{k-1}, 0) + b \min(X_{k-1}, 0) + \varepsilon_k,$$

where $a$ and $b$ are real parameters. For processes of the form (2.10), condition (2.5) can be simplified. Let $\varphi$ be the characteristic function of $\varepsilon_1$.

**Theorem 2.2.** *Let $0 < \alpha \leq 2$. Assume (2.10),*

$$\int_{\mathbb{R}} |\varphi(t)|^2 (1+t^2) |t|^\alpha dt < \infty \tag{2.12}$$

*and*

$$\sum_{k=0}^{\infty} \|X_k - X_k^*\|_\alpha^{\alpha/2} < \infty \tag{2.13}$$

*Then (2.5) is satisfied.*

Condition (2.12) does not appear to be overly restrictive. It is satisfied if $|\varphi(t)| = O(|t|^{-\eta})$ as $|t| \to \infty$, where $\eta > (3+\alpha)/2$. It is also satisfied for symmetric-$\alpha$-stable distributions, an important class of distributions with heavy tails. Let $\varepsilon_k$ have standard symmetric $\alpha$ stable distribution with index $0 < \iota \leq 2$. Then its characteristic function $\varphi(t) = \exp(-|t|^\iota)$ and (2.12) trivially holds.

We now discuss Condition (2.13). Recall $X_k^* = g(\xi_k^*)$. Note that $X_k^*$ and $X_k$ are identically distributed and $X_k^*$ is a coupled version of $X_k$ with $\varepsilon_0$ replaced by $\varepsilon_0'$. If we view (1.1) as a physical system with $\xi_k = (\ldots, \varepsilon_{k-1}, \varepsilon_k)$ being the input, $g$ being a filter or tranform and $X_i$ being the output, then the quantity $\|X_k - X_k^*\|_\alpha$ measures the degree of dependence of $g(\ldots, \varepsilon_{k-1}, \varepsilon_k)$ on $\varepsilon_0$. In [23] it is called the *physical* or *functional dependence measure*. With this input/output viewpoint, the condition (2.13) means that the cumulative impact of $\varepsilon_0$ is finite, and hence suggesting short-range dependence. In many applications it is easily verifiable since it is directly related to the data-generating mechanism and since the calculation of $\|X_k - X_k^*\|_\alpha$ is generally easy [23]. In the special case of linear process $X_k = \sum_{j=0}^{\infty} a_j \varepsilon_{k-j}$ with $\varepsilon_k \in \mathcal{L}^\alpha$ and $\alpha = 2$, then $\|X_k - X_k^*\|_\alpha = |a_k| \|\varepsilon_0 - \varepsilon_0'\|_\alpha$ and (2.13) is reduced to $\sum_{k=0}^{\infty} |a_k| < \infty$, which is a classical condition for linear processes to be short-range dependent. It is well-known that, if the latter condition is barely violated, then one enters the territory of long-range dependence. Consequently both the normalization and the bound in (2.7) will be different; see [8, 21].

For the nonlinear time series (2.11), assume that $\varepsilon_k \in \mathcal{L}^\alpha$ and $\rho = \sup_x |m'(x)| < 1$. Then (2.11) has a stationary distribution and $\|X_k - X_k^*\|_\alpha = O(\rho^k)$ (see [24]). Hence (2.13) holds.

### 3. Proofs

**Lemma 3.1.** *Let $H$ be a differential function on $\mathbb{R}$. Then for any $\lambda > 0$,*

$$\sup_{x \in \mathbb{R}} H^2(x) \leq \lambda \int_{\mathbb{R}} H^2(x) dx + \lambda^{-1} \int_{\mathbb{R}} [H'(x)]^2 dx. \tag{3.1}$$

*Proof.* By the arithmetic mean geometric inequality inequality, for all $x, y \in \mathbb{R}$,

$$H^2(x) \leq H^2(y) + \left|\int_x^y 2H(t)H'(t)dt\right|$$

(3.2)
$$\leq H^2(y) + \lambda \int_{\mathbb{R}} H^2(x)dx + \lambda^{-1} \int_{\mathbb{R}} [H'(x)]^2 dx.$$

If $\inf_{y \in \mathbb{R}} |H(y)| > 0$, then $\int_{\mathbb{R}} H^2(x)dx = \infty$ and (3.1) holds. If on the other hand $\inf_{y \in \mathbb{R}} |H(y)| = 0$, let $(y_n)_{n \in \mathbb{N}}$ be a sequence such that $H(y_n) \to 0$. So (3.2) entails (3.1). □

Lemma 3.1 is a special case of the Kolmogorov-type inequalities [18]. The result in the latter paper asserts that $\sup_{x \in \mathbb{R}} H^4(x) \leq \int_{\mathbb{R}} H^2(x)dx \times \int_{\mathbb{R}} H'(x)^2 dx$. For the sake of completeness, we decide to state Lemma 3.1 with a simple proof here.

Recall that $F_k(\cdot|\xi_0)$ is the conditional distribution function of $X_k$ given $\xi_0$ (cf (2.1) and (2.2)). Introduce the conditional empirical distribution function

$$F_n^*(x) = \frac{1}{n} \sum_{i=1}^n \mathbb{E}(\mathbf{1}_{X_i \leq x}|\xi_{i-1}) = \frac{1}{n} \sum_{i=1}^n F_1(x|\xi_{i-1}).$$

Write

(3.3)
$$G_n(x) = G_n^\diamond(x) + G_n^*(x),$$

where

(3.4)    $G_n^\diamond(x) = \sqrt{n}[F_n(x) - F_n^*(x)]$ and $G_n^*(x) = \sqrt{n}[F_n^*(x) - F(x)].$

Then

(3.5)
$$\sqrt{n}G_n^\diamond(x) = \sum_{i=1}^n d_i(x)$$

is a martingale with respect to the filtration $\sigma(\xi_n)$ and the increments $d_i(x) = \mathbf{1}_{X_i \leq x} - \mathbb{E}(\mathbf{1}_{X_i \leq x}|\xi_{i-1})$ are stationary, ergodic and bounded. On the other hand, if the conditional density $f_1(\cdot|\xi_i)$ exists, then $G_n^*$ is differentiable. The latter differentiability property is quite useful.

**Lemma 3.2.** *Recall (2.6) for $\iota(n)$. Let $g_n^*(x) = dG_n^*(x)/dx$. Assume (2.5). Then*

(3.6)
$$\sup_x |g_n^*(x)| = o_{a.s.}[\iota(n)].$$

*Proof.* Let $k \geq 1$. Recall that $f_1(x|\xi_{k-1})$ is the one-step-ahead conditional density of $X_k$ at $x$ given $\xi_{k-1}$. By (2.3), we have

$$\mathcal{P}_0 \varphi_1(t|\xi_{k-1}) = \int_{\mathbb{R}} e^{\sqrt{-1}xt} \mathcal{P}_0 f_1(x|\xi_{k-1}) dt.$$

By Parseval's identity, we have

$$\int_{\mathbb{R}} |\mathcal{P}_0 \varphi_1(t|\xi_{k-1})|^2 dt = \frac{1}{2\pi} \int_{\mathbb{R}} |\mathcal{P}_0 f_1(x|\xi_{k-1})|^2 dx$$

and

$$\int_{\mathbb{R}} |\mathcal{P}_0 \varphi_1(t|\xi_{k-1})|^2 t^2 dt = \frac{1}{2\pi} \int_{\mathbb{R}} |\mathcal{P}_0 f_1'(x|\xi_{k-1})|^2 dx.$$

Let
$$\alpha_k = \int_{\mathbb{R}} \|\mathcal{P}_0 f_1(x|\xi_{k-1})\|^2 dx$$

and
$$\beta_k = \int_{\mathbb{R}} \|\mathcal{P}_0 f_1'(x|\xi_{k-1})\|^2 dx.$$

By (2.8) and Jensen's inequality, $\|\mathcal{P}_0 \varphi_1(\theta|\xi_{k-1})\| \leq \|\varphi_k(\theta|\xi_0) - \varphi_k(\theta|\xi_0^*)\|$. So (2.5) implies that

(3.7) $$\sum_{k=1}^{\infty} \sqrt{\alpha_k + \beta_k} < \infty.$$

Let $\Lambda = \sum_{k=1}^{\infty} \sqrt{\alpha_k}$ and for $k \geq 0$
$$H_k(x) = \sum_{i=1}^{k} [f_1(x|\xi_{i-1}) - f(x)].$$

Then $H_k(x) = \sqrt{k} g_k(x)$. Note that for fixed $l \in \mathbb{N}$, $\mathcal{P}_{i-l} f_1(x|\xi_{i-1})$, $i = 1, 2, \ldots$, are stationary martingale differences with respect to the filtration $\sigma(\xi_{i-l})$. By the Cauchy-Schwarz inequality and Doob's maximal inequality, since $H_k(x) = \sum_{l=1}^{\infty} \sum_{i=1}^{k} \mathcal{P}_{i-l} f_1(x|\xi_{i-1})$,

$$\mathbb{E} \int_{\mathbb{R}} \max_{k \leq n} H_k^2(x) dx \leq \mathbb{E} \int_{\mathbb{R}} \sum_{l=1}^{\infty} \frac{\max_{k \leq n} |\sum_{i=1}^{n} \mathcal{P}_{i-l} f_1(x|\xi_{i-1})|^2}{\sqrt{\alpha_l}} \Lambda dx$$
$$\leq \int_{\mathbb{R}} \sum_{l=1}^{\infty} \frac{4n \|\mathcal{P}_0 f_1(x|\xi_{l-1})\|^2}{\sqrt{\alpha_l}} \Lambda dx = 4n\Lambda^2 = O(n).$$

Similarly, since $\sum_{k=1}^{\infty} \sqrt{\beta_k} < \infty$,
$$\mathbb{E} \int_{\mathbb{R}} \max_{k \leq n} |H_k'(x)|^2 dx = O(n).$$

By Lemma 3.1 with $\lambda = 1$, we have
$$\sum_{d=1}^{\infty} \frac{\mathbb{E}[\max_{k \leq 2^d} \sup_{x \in \mathbb{R}} |H_k(x)|^2]}{2^d \iota^2(2^d)}$$
$$\leq \sum_{d=1}^{\infty} \frac{\mathbb{E}[\max_{k \leq 2^d} \int_{\mathbb{R}} |H_k(x)|^2 + |H_k'(x)|^2 dx]}{2^d \iota^2(2^d)}$$
$$\leq \sum_{d=1}^{\infty} \frac{\mathbb{E} \int_{\mathbb{R}} \max_{k \leq 2^d} |H_k(x)|^2 + \max_{k \leq 2^d} |H_k'(x)|^2 dx}{2^d \iota^2(2^d)}$$
$$= \sum_{d=1}^{\infty} \frac{O(2^d)}{2^d \iota^2(2^d)} < \infty.$$

By the Borel-Cantelli lemma, $\sup_x \max_{k \leq 2^d} |H_k(x)| = o_{\text{a.s.}}[2^{d/2} \iota(2^d)]$ as $d \to \infty$. For any $n \geq 2$ there is a $d \in \mathbb{N}$ such that $2^{d-1} < n \leq 2^d$. Note that $\max_{k \leq n} |H_k(x)| \leq \max_{k \leq 2^d} |H_k(x)|$ and $\iota(n)$ is slowly varying. So (3.6) follows. □

**Lemma 3.3.** *Assume* $\log n = O(nb_n)$ *and* $X_k \in \mathcal{L}^\alpha$ *for some* $\alpha > 0$. *Then for any* $\tau > 2$, *there exists* $C = C_\tau > 0$ *such that*

$$\mathbb{P}\left[\sup_{|x-y|\leq b_n} |G_n^\diamond(x) - G_n^\diamond(y)| > C\sqrt{b_n \log n}\right] = O(n^{-\tau}). \tag{3.8}$$

*Proof.* We can adopt the argument of Lemma 5 in [22]. Let $x_0 = n^{(3+\tau)/\alpha}$. Then for any $C > 0$, by Markov's inequality,

$$\begin{aligned}
&\mathbb{P}[n\{F_n(-x_0) + 1 - F_n(x_0)\} > c\sqrt{b_n \log n}] \\
&\leq \frac{n\mathbb{E}\{F_n(-x_0) + 1 - F_n(x_0)\}}{C\sqrt{b_n \log n}} \\
&\leq \frac{nx_0^{-\alpha}\mathbb{E}(|X_0|^\alpha)}{C\sqrt{b_n \log n}} = O(n^{-\tau}).
\end{aligned} \tag{3.9}$$

It is easily seen that the preceding inequality also holds if $F_n$ is replaced by $F_n^*$.

Recall (3.5) for $G_n^\diamond(x)$ and $d_i(x) = \mathbf{1}_{X_i \leq x} - \mathbb{E}(\mathbf{1}_{X_i \leq x}|\xi_{i-1})$. Let $x \leq y \leq x + b_n$. By (2.4),

$$\sum_{i=1}^n \mathbb{E}[(d_i(y) - d_i(x))^2 | \xi_{i-1}] \leq \sum_{i=1}^n \mathbb{E}(\mathbf{1}_{x \leq X_i \leq y}|\xi_{i-1}) \leq nb_n c_0.$$

By Freedman's inequality in [6], if $|x - y| \leq b_n$, we have

$$\mathbb{P}\left[\sqrt{n}|G_n^\diamond(x) - G_n^\diamond(y)| > C\sqrt{nb_n \log n}\right] \leq 2 \exp\left[\frac{-C^2 nb_n \log n}{C\sqrt{nb_n \log n} + nb_n c_0}\right].$$

Let $\Theta_n = \{-x_0 + k/n^3 : k = 0, \ldots, \lfloor 2x_0 n^3 \rfloor\}$. Since $\log n = O(nb_n)$, it is easily seen that there exists a $C = C_\tau$ such that

$$\begin{aligned}
&\mathbb{P}\left[\sup_{x,y \in \Theta_n, |x-y| \leq b_n} \sqrt{n}|G_n^\diamond(x) - G_n^\diamond(y)| > C\sqrt{nb_n \log n}\right] \\
&= O(x_0^2 n^2) \exp\left[\frac{-C^2 nb_n \log n}{C\sqrt{nb_n \log n} + nb_n c_0}\right] = O(n^{-\tau}).
\end{aligned} \tag{3.10}$$

For every $x \in [-x_0, x_0]$, there exists a $\theta \in \Theta_n$ such that $\theta < x \leq \theta + 1/n^3$. So (3.10) implies that

$$\begin{aligned}
&\mathbb{P}\left[\sup_{|x| \leq x_0, |y| \leq x_0, |x-y| \leq b_n} \sqrt{n}|G_n^\diamond(x) - G_n^\diamond(y)| > (C+1)\sqrt{nb_n \log n}\right] \\
&= O(n^{-\tau})
\end{aligned} \tag{3.11}$$

in view of the monotonicity of $F_n(\cdot)$ and the fact that, if $\theta \leq \phi \leq \theta + 1/n^3$,

$$\sum_{i=1}^n \mathbb{E}(\mathbf{1}_{\theta \leq X_i \leq \phi}|\xi_{i-1}) \leq |\phi - \theta|nc_0 = c_0/n^2 = o(\sqrt{nb_n \log n}).$$

Combining (3.9) and (3.11), we have (3.8). □

*Proof of Theorem 2.1.* By (3.3), it easily follows from Lemmas 3.2 and 3.3. □

*Proof of Theorem 2.2.* By Lemma 3.1 and Parseval's identity, (2.12) implies

$$\sup_x f^2(x) \leq \int_{\mathbb{R}} \{f^2(x) + [f'(x)]^2\}dx = \frac{1}{2\pi}\int_{\mathbb{R}} |\varphi(t)|^2(1+t^2)dt < \infty.$$

Since $\varepsilon_k$ and $Y_{k-1}$ are independent,

$$\mathbb{E}[e^{\sqrt{-1}tX_k}|\xi_{k-1}] = e^{\sqrt{-1}tY_{k-1}}\varphi(t).$$

Note that

$$\begin{aligned}
&\|e^{\sqrt{-1}tY_{k-1}}\varphi(t) - e^{\sqrt{-1}tY^*_{k-1}}\varphi(t)\|^2 \\
&= |\varphi(t)|^2 \|e^{\sqrt{-1}tY_{k-1}} - e^{\sqrt{-1}tY^*_{k-1}}\|^2 \\
&\leq 4|\varphi(t)|^2 \|\min(1, |tY_{k-1} - tY^*_{k-1}|)\|^2 \\
&\leq 4|\varphi(t)|^2 \mathbb{E}(|tY_{k-1} - tY^*_{k-1}|^\alpha).
\end{aligned}$$

Hence (2.5) follows from (2.13). □

## 4. Conclusion and open problems

Let $X_i$ be iid standard uniform random variables. Stute [16] obtained the following interesting result. Assume that $b_n \to 0$ is a sequence of positive numbers such that

(4.1) $$\log n = o(nb_n) \text{ and } \log\log n = o(\log b_n^{-1}).$$

Then the convergence result holds:

(4.2) $$\lim_{n\to\infty} \frac{\Delta_n(b_n)}{\sqrt{b_n \log b_n^{-1}}} = \sqrt{2} \text{ almost surely.}$$

If there exists $\eta > 0$ such that $b_n + (nb_n)^{-1} = O(n^{-\eta})$, then the bound in (2.7) becomes $\sqrt{b_n \log n}$, which has the same order of magnitude as $\sqrt{b_n \log b_n^{-1}}$, the bound asserted by (4.2). It is unclear whether there exists an almost sure limit for $\Delta_n(b_n)/\sqrt{b_n \log b_n^{-1}}$ if the dependence among observations is allowed. Mason et al [11] considered almost sure limit for $\Delta_n(b_n)/\sqrt{b_n \log b_n^{-1}}$ when (4.1) is violated. Deheuvels and Mason [3] (see also [2]) proved functional laws of the iterated logarithm for the increments of empirical processes. Local empirical processes in high dimensions have been studied in [5, 7, 17]. It is an open problem whether similar results hold for stationary causal processes. We expect that our decomposition (3.4) will be useful in establishing comparable results.

### Acknowledgements

I am grateful for the very helpful suggestions from the reviewer and the Editor.